\documentclass[11pt]{article}

\usepackage{amsmath,amssymb}

\newtheorem{theorem}{Theorem}
\newtheorem{lemma}[theorem]{Lemma}
\newtheorem{prop}[theorem]{Proposition}
\newtheorem{corr}[theorem]{Corollary}

\newcommand{\X}{{\mathbb{X}_{\Sigma}}}
\newcommand{\R}{{\mathbb{R}}}

\newcommand{\Z}{{\mathbb{Z}}}
\newcommand{\s}[1]{\widehat{#1}}
\renewcommand{\k}{k}

\begin{document}

\title{BGG correspondence for toric complete intersections}
\author{Vladimir Baranovsky}
\date{April 14, 2007}
\maketitle

\begin{abstract}
We prove a BGG-type correspondence describing coherent sheaves
on complete intersections in toric varieties, and  a
similar assertion for the stable categories of related complete
intersection singularities.
\end{abstract}

\section{Introduction}

This paper is a continuation of the earlier article on complete 
intersections in projective spaces, cf. \cite{Ba}. We consider
here the case of a complete intersection $Y$ in a 
toric variety $X_\Sigma$ over a field $k$ 
of characteristic zero. In the case when $X_\Sigma$ has singularities,
we actually study the corresponding \textit{stacks} $\mathbb{Y}
\subset \X$ (this point of view is also used, for instance,
when toric complete intersections are considered in Mirror 
Symmetry). Our goal is to give an alternative 
description for the category of sheaves on such a $\mathbb{Y}$
in the spirit of the one given by 
Bernstein-Gelfand-Gelfand in \cite{BGG} for projective spaces
and by Kapranov in \cite{Ka} for intersections of projective
quadrics. The general approach is modeled on the Koszul duality
of Beilinson-Ginzburg-Schechtman, cf. \cite{BGS}, but in 
our case we deal with the higher products on 
the ``Koszul dual" which arise from the fact that the original
algebra had non-quadratic relations. 

Now we describe the contents in more detail.
In the above setting $\X$ has a ``homogenous coordinate ring" $S$
isomorphic to a polynomial algebra graded by a finitely generated abelian
group $A$, cf. \cite{C1}. If $W_1, \ldots, W_m$ are the defining equations of
$\mathbb{Y}$ and $J$ is the ideal of $S$ generated by these
equations, then the category $Coh(\mathbb{Y})$ of coherent
sheaves on $\mathbb{Y}$ is obtained from the category 
of finitely generated 
$A$-graded modules over $S_W = S/J$ by passing to a certain 
categorical quotient, see Section 4 for details.

We first study $A$-graded modules over $S_W$. In Section 2
we use the polynomials $W_1, \ldots, W_m$  to define, more or 
less tautologically, an $L_\infty$-algebra $L$.
We further construct an $A_\infty$-algebra $E_W$ which 
should be viewed as the ``universal enveloping" of $L$. 
When $W_1, \ldots, W_m$ have no linear terms (which one 
can alsways assume replacing $S$ by a quotient polynomial
algebra),
$E_W$ has zero differential. In the case when all $W_j$
are quadratic $E_W$ becomes the associative graded Clifford
algebra considered by Kapranov in \cite{Ka}. The proof 
proceeds differently from \cite{Ba} since we do not assume
that $W_1, \ldots, W_m$ are homogeneous with respect to the 
usual grading on $S$ (which is necessary for toric applications).
Ideally, one would like to characterize $E_W$ as the unique
``homotopy bialgebra" of some special sort, 
such that the restriction of  $A_\infty$-products to 
$L \subset E_W$ is given by the homogeneous components of 
$W_1, \ldots, W_m$. However, we leave the task of writing the
agreement conditions between the $A_\infty$-products and the
natural coproduct on $E_W$, to a forthcoming paper.

In Section 3 we prove, see Theorem \ref{coequiv}, 
an equivalence between $A_\infty$-modules 
over $E_W$ and $L_\infty$-modules over $L$ (the latter are 
viewed as modules over the standard Cartan-Eilenberg-Chevalley 
coalgebra $C$  of $L$). The proof follows the general 
formalizm developed in \cite{Le} which we expand slightly 
to the $A_\infty$-case. When $S_W$ is graded by $A$ as above
and all graded components are finite-dimensional, we obtain
an equivalence between $A$-graded modules over $S_W$ and 
$E_W$, see Theorem \ref{bgg}. This result is applied in 
Section 4 to the derived category of sheaves on a toric 
complete intersection  $\mathbb{Y}$ and to the stable category 
of the affine complete intersection defined by $W_1, \ldots, W_m$. 
When  $\mathbb{Y}$ has trivial canonical class (with an 
additional technical assumption always satisfied 
for intersections in weighted projective spaces)
an easy application of a result due to Orlov, cf.
\cite{O2}, gives an 
alternative description of the derived category of $\mathbb{Y}$,
cf. Corollary 12.

\bigskip
\noindent
\textbf{Acknowledgements.} This work was supported by the Sloan
Research Fellowship.

\section{A universal enveloping algebra}

\subsection{Differential operators and corrected partial derivatives}

Fix a finite dimensional vector space $V$ over $k$. The symmetric
algebras $Sym^\bullet(V)$ and $Sym^\bullet(V^*)$ may be viewed
as algebras of differential operators (with constant coefficients) over each 
other. For any $f \in Sym^\bullet(V)$ let $\partial_f$ be the corresponding
operator on $Sym(V^*)$, and similarly for $g \in Sym^\bullet(V^*)$.
There is a pairing 
$\langle \cdot, \cdot\rangle: Sym^\bullet(V) \times 
Sym^\bullet(V^*) \to k$ given by
$$
\langle f, g\rangle := \partial_f(g)(0) = \partial_g (f)(0)
$$
With respect to this pairing, $\partial_g$ is adjoint
to multiplication by $g$ on $Sym^\bullet(V^*)$. 

We will also need ``corrected partial derivatives"': for $v \in V$ let 
$\widehat{\partial}_v$ be the operator which sends $g \in Sym^k(V^*)$ to 
$\frac{1}{k} \partial_v(g)$ for $k \geq 1$ and satisfies 
$\widehat{\partial}_v (1) = 0$.

For a vector space $U$ we view $Sym^\bullet(V)$ as differential 
operators on $Sym^\bullet(V^*) \otimes U$ extending derivatives
(usual or ``corrected") by linearity in the 
second factor.

\subsection{Koszul complex and an $L_\infty$-algebra.}

Choose and fix a regular sequence  $W_1, \ldots, W_m
\in Sym^{\geq 1}(V^*)$.
Introducing new variables $z_1, \ldots, z_m$ which span a vector space
$U$ we can encode the above sequence in a single ``total potential"
$$
W = W_1 z_1 + \ldots + W_m z_m \in Sym^\bullet(V^*) \otimes U
$$
Unlike in 
\cite{Ba}, we do not make the assumption that $W_j$ are homogeneous.

Due to the regularity, the quotient $S_W = Sym^\bullet(V^*)/J$ 
by the ideal generated $J$ by $W_j$, $j = 1, \ldots, m$, admits 
a Koszul resolution $B = Sym^\bullet(V^*) \otimes \Lambda^\bullet(U^*)$
where the differential $\delta_B$ is given by $W$, if we agree that $z_j$
act on $\Lambda^\bullet(U^*)$ by contraction and $W_j$ on $Sym^\bullet(V^*)$
by multiplication.

The differential $\delta_C$ of the dual coalgebra 
$C = Sym^\bullet(V) \otimes \Lambda(U)$ is also given by $W$ but 
now we think of $W_j$ as differential operators and $z_j$ act by 
multiplication (in the natural algebra structure on $C$). The
assumption that the sequence $W_1, \ldots, W_m$ is regular
will not be needed in this Section.

Introduce an $L_\infty$-algebra $L = \big\{0 \to V \to  U\big\}$,
cf. \cite{LM}, placing $V$  in 
homological degree 1, $U$ in homological degree 2 and defining 
the $L_\infty$-operations as follows. We set
$$
l_k(v_1, \ldots, v_k) := \partial_{v_1 \ldots v_k} (W)(0) = k! \;
\widehat{\partial}_{v_1} \ldots \widehat{\partial}_{v_k}(W)(0)
$$
whenever all arguments 
$v_1 \ldots, v_k$ are in $V$, and let $l_k = 0$ otherwise. 
The $L_\infty$-identities for $L$ will hold trivially, since 
every double composition involved in them vanishes. Note
that the coproduct of $C$ is independent of $W$, but its
differential contains full information about it. Also,
 $L$ is isomorphic as a vector space to the space of 
primitive elements in $C$. In fact, one has the following lemma 
which is immediate from definitions
\begin{lemma}\label{coalgebra}
The coalgebra $C$ is isomorphic to the  cocommutative coalgebra $C(L)$ of
$L$, cf. \cite{LM}. 
\end{lemma}

\subsection{A standard resolution of $L$}

We now describe a resolution $L \to \mathcal{L}$ 
 in which the bracket does not depend on $W$. 
Let $\mathcal{L}$ be the graded vector space with 
$\mathcal{L}^1 = \Big[Sym^{\geq 1} (V^*) \otimes U\Big] 
\oplus V$ in homological degree 1 and
$\mathcal{L}^2 = Sym^\bullet(V^*) \otimes U$ in homological 
degree 2. Define the differential 
\begin{equation}\label{differential}
\delta_\mathcal{L} (f \oplus v) = \widehat{\partial}_v(W) - [f]
\end{equation}
where  for $f \in Sym^{\geq 1}(V^*) \otimes U \subset \mathcal{L}^1$ we 
denote by $[f]$ its copy in $\mathcal{L}^2$. The bracket
$\{\;,\;\}: \mathcal{L}^1 \times \mathcal{L}^1 \to \mathcal{L}^2$ 
is defined by 
\begin{equation}\label{bracket}
\{f_1 \oplus v_1, f_2 \oplus v_2\} = \widehat{\partial}_{v_1}(f_2) 
+  \widehat{\partial}_{v_2}(f_1)
\end{equation}
For any $w \in Sym^\bullet(V^*) \otimes U$ denote
by $\overline{w} \in Sym^{\geq 1}(V^*)\otimes U$ its image with 
respect to 
the natural projection which has $k \otimes U$ as its kernel. 
Then the morphism of complexes $G_1: L \to \mathcal{L}$
$$
G_1(v) = \overline{\partial_v (W)} + v \in 
\Big[Sym^{\geq 1} (V^*) \otimes U\Big]  \oplus V, 
\quad G_1(u) = u \in U \subset \mathcal{L}^2;
\qquad v \in V = L^1, u \in U = L^2
$$
is a quasi-isomorphism, but not a morphism of DG Lie algebras. 
However, introducing morphisms
$$
G_k: L^k \to \mathcal{L}, \qquad (v_1, \ldots, v_k) 
\mapsto k!\;\overline{\widehat{\partial}_{v_1} \ldots 
\widehat{\partial}_{v_k}(W)}
$$
whenever all $v_i$ are in $V$, and zero otherwise; we extend 
$G_1$ to an $L_\infty$-morphism $\{G_k\}_{k \geq 1}$, 
cf. \cite{LMa}. The $L_\infty$-morphism condition of 
\textit{loc. cit.} in our case reduces to 
\begin{equation}\label{identity}
G_1(l_k(v_1, \ldots, v_k)) + \delta_{\mathcal{L}}(G_k(v_1, \ldots, v_k))
= \sum_{i=1}^k \{G_{k-1}(v_1, \ldots, \widehat{v_i}, \ldots, v_k), G_1(v_i)\}
\end{equation}
when $k \geq 3$; while for $k = 2$ one has
$$
G_1(l_2(v_1, v_2)) + \delta_{\mathcal{L}}(G_2(v_1, v_2)) = \{G_1(v_1), G_1(v_2)\}
$$
We note here that it is precisely \eqref{identity} why we use
``corrected partial derivatives" in the definitions of $l_k$ and $G_k$.

The individual maps $G_k$, $k \geq 1$ can be organized into a
single map $G_\infty: C\to Sym^{\geq 1}(V) \to \mathcal{L}$. 
Since $C$ is a 
cocommutative coalgebra, by Lemma 22.1 in \cite{FHT} 
there is a unique comultiplicative extension
$\tau: C \to Sym^\bullet_c(\mathcal{L})$ into the symmetric coalgebra of 
$\mathcal{L}$. We further use Poincare-Birkhoff-Witt 
to identify $Sym^\bullet_c(\mathcal{L})$ 
with the universal enveloping $U(\mathcal{L})$ 
of $\mathcal{L}$ (as DG coalgebras).

The following lemma deals the multiplicative 
behavior of $\tau$ with respect to the standard universal
enveloping product $m_2$ in of $U(\mathcal{L})$ and the 
product in the reduced cobar construction $\Omega(C)$. See 
e.g. \cite{FHT} and \cite{Ka2} regarding the 
definitions and properties of the cobar construction.

\begin{lemma}\label{quotient}
The unique comultiplicative extension 
$\tau: C \to Sym^\bullet_c(\mathcal{L}) 
\simeq U(\mathcal{L})$ satisfies the twisted cochain condition
$$
\tau \circ \delta_C + \delta_{U(\mathcal{L})} \circ \tau 
+ m_2 \circ \tau^{\otimes 2} \circ \Delta =0.
$$
Its own canonical multiplicative extension 
$\Omega(\tau): \Omega(C) \to U(\mathcal{L})$ 
is a quasi-isomorphsim of DG algebras.
\end{lemma}
\textit{Proof.} By \eqref{identity} above the map 
$G_\infty: C \to \mathcal{L}$
extends to a morphism of DG-coalgebras $C \to C(\mathcal{L})$,
where $C(\cdot)$ stands for the Cartan-Eilenberg-Chevalley 
coalgebra of a DG Lie algebra, cf. \cite{FHT}. 
It is an easy computation that the composition of natural
maps $C(\mathcal{L}) \to \mathcal{L} \to U(\mathcal{L})$
does satisfy the twisted cochain condition. 
Since $\tau: C \to U(\mathcal{L})$ factors 
as $C \to C(\mathcal{L}) \to U(\mathcal{L})$, 
the first assertion follows.

For the second assertion note that $\Omega(C) \to U(\mathcal{L})$
 commutes with differentials due to the twisted cocycle property 
of $\tau$. 

For any DG coalgebra $C'$ let $\mathbb{L}(C')$
be Quillen's free DG Lie algebra of $C'$, cf. Section 22(e) of
\cite{FHT}. Then $\Omega(C') \simeq U\mathbb{L}(C')$ by the 
universal properties of the three objects involved.
Now decompose $\Omega(\tau)$ as 
$\Omega(C) \to \Omega(C(\mathcal{L}))
= U \mathbb{L}(C(\mathcal{L})) \to U(\mathcal{L})$. The 
first arrow is a quasi-isomorphism because the $L_\infty$ map 
$G_\infty: C= C(L) \to \mathcal{L}$ extends to
a quasi-isomorphism of DG coalgebras (this follows from the 
fact that $G_1$ is a quasi-isomorphism of complexes). 
The second arrow is a quasi-isomorphism since it is induced by 
a quasi-isomorphism of DG Lie algebras
$\mathbb{L}(C(\mathcal{L})) \to \mathcal{L}$, cf. 
Theorem 22.9 in \cite{FHT}.
$\square$.

\subsection{The universal enveloping $A_\infty$-algebra $E_W$}

By  Theorem 22.9 in \cite{FHT}
for a DG Lie algebra $L'$ one has a quasi-isomorphism of DG 
algebras $\Omega(C(L')) \to U(L')$. We want to use this fact to 
define an $A_\infty$-structure on the symmetric coalgebra 
$Sym_c(L)$ which should be viewed as the ``universal enveloping" 
algebra of $L$. An ideal strategy would be as follows: first 
replace $W$ in the definition of $L$ with the potential $W^{(2)}$ 
obtained from $W$ by erasing the terms of degree $\geq 3$ in the 
usual homogeneous grading of $Sym^\bullet(V^*)$. In other words,
we forget all higher brackets on $L$ which in our case leads
to a DG Lie algebra $L_2$ and a quasi-isomorphism 
$\Omega(C(L_2)) \to U(L_2)$. Bringing back the degree $\geq 3$
terms of $W$ amounts to perturbing the differential on 
$\Omega(C(L_2))$ and the ``sum over binary trees" formula
of \cite{KS} tells us that this perturbation induces 
an $A_\infty$-structure on $U(L_2)$. 

However, this formula involves an explicit contracting 
homotopy on $\Omega(C(L_2))$ which we are not able to write
down at the moment. Therefore we replace $\Omega(C)$ by 
a smaller DG algebra $U(\mathcal{L})$ which is quasi-isomorphic
to it by Lemma \ref{quotient}. Moreover, we do not
apply the ``sum over binary trees" formula but rather the 
results of \cite{GLS} which, in a sense, stand behind it. 
In more detail: replacing $W$ by $W^{(2)}$ gives a DG Lie
algebra $\mathcal{L}_2$ and a quasi-isomorphism of DG Lie
algebras $G_1: L_2 \to \mathcal{L}_2$ (all higher $G_k$ vanish 
in for $W^{(2)}$). We also denote by 
$\mathcal{L}_1$ and $L_1$  the same objects viewed 
 as complexes with trivial Lie bracket. 
First we construct a canonical contracting homotopy on 
$U(\mathcal{L}_1)$
and then take into account the Lie brackets and use
\cite{GLS} to define an $A_\infty$-map $\{F_k\}_{k \geq 1}$
of associative algebras $U(\mathcal{L}_2) \to U(L_2)$ and
a system of higher homotopies $\{H_k\}_{k \geq 1}$ 
on $U(\mathcal{L}_2)$. Finally, we replace $W^{(2)}$ by $W$
and then the constructed system of  homotopies gives
an $A_\infty$-structure on $E_W$. 

The advantage of this approach, which is more complicated than
fixing a non-canonical homotopy on $U(\mathcal{L}_2)$, is 
that the resulting $A_\infty$ structure on $U(L_2)$ only
depends on the resolution $\mathcal{L}$ and, in addition,
it has some compatibility with the coproduct (see the remark
at the end of this section).

\bigskip
\noindent
So replace $W$ by $W^{(2)}$ as above and consider
the complexes  $L_1, \mathcal{L}_1$. 
The map $G_1: L_1 \to \mathcal{L}_1$
of the previous section admits a left inverse 
$F: \mathcal{L}_1 \to L_1$ which projects 
$\Big[Sym^{\geq 1} (V^*) \otimes U\Big]  \oplus V 
= (\mathcal{L}_1)^1$
onto $V = (L_1)^1$ in an obvious way (the superscripts
denote homological grading), and sends 
$ Sym^\bullet(V^*)\otimes U  = (\mathcal{L}_1)^2$ to 
$U = (L_1)^2$ by evaluating the constant term. 
Define a homotopy $H: \mathcal{L}^2_1 \to \mathcal{L}^1_1$
by sending $w$ to $\{\overline{w}\}$  (we use braces to  
emphasize that an even element $w$ was converted into an odd 
element).  The ``side conditions"
$$
HG_1 = 0, \quad HH = 0, \quad FH=0
$$  
follow immediately from the definitions.

Now we consider the symmetric DG bialgebras (in the graded 
sense) $Sym^\bullet(\mathcal{L}_1)$ and 
$Sym^\bullet(L_1)$ and the natural extensions of $F$ and $G_1$ 
given by multiplicative and comultiplicative maps $F_{sym}:  Sym^\bullet(\mathcal{L}_1)
\to Sym^\bullet(L_1)$, and $G_{sym}:  Sym^\bullet(L_1)
\to Sym^\bullet(\mathcal{L}_1)$. To define a homotopy $H_{sym}$ we 
first set  $\mathcal{S} = Sym^{\geq 1}(V^* \otimes U) $ and denote
by  $\{S\}$, $[S]$ its copies in $\mathcal{L}^1_1$ and $\mathcal{L}^2_1$, 
respectively. Since $\mathcal{L}_1 = G_1(L_1) \oplus \big(\{S\} \to [S]\big)$
as complexes, we have an isomorphism of DG bialgebras
$$
Sym^\bullet(\mathcal{L}_1)  \simeq Sym^\bullet(\{S\} \to [S]) \otimes
Sym^\bullet(L_1).
$$
The graded symmetric bialgebra $Sym^\bullet(\{S\} \to [S]) \simeq \Lambda^\bullet(S)  \otimes Sym^\bullet(S)$
has standard Koszul differential, and therefore a standard 
homotopy 
\begin{equation}\label{homotopy}
H_{sym}\Big(\{f_1\}  \ldots  \{f_m\}
[g_1]  \ldots  \;[g_k]\Big) 
= \frac{-1}{k+m} \sum_{t = 1}^k \{f_1\}  \ldots  \{f_m\}  \{g_t\}
[g_1]  \ldots \widehat{[g_t]}  \ldots  [g_k].
\end{equation}
This we extend to $Sym^\bullet(\mathcal{L}_1)$ as 
$H_{sym} \otimes 1$ denoting the extension again by $H_{sym}$. 

The contraction $(F_{sym}, G_{sym}, H_{sym})$ induces a 
similar contraction $(F_B', G_B', H_B')$ on the
reduced bar constructions (see Section 19 of \cite{FHT}
and \cite{Ka2} for definitions and properties). Here
$F_B', G_B'$ are defined in an obvious way and 
$$
H_B'|_{(Sym^{\geq 1}(\mathcal{L}_1))^{\otimes k}} = \sum_{s = 0}^{k-1} 
1^{\otimes s} \otimes H_{sym} \otimes (G_{sym} F_{sym})^{\otimes(k-s-1)}
$$
Then  $F_B'$ and $G_B'$ are maps of DG coalgebras and $H_B'$ is a
coalgebra homotopy:
\begin{equation}\label{ch}
\Delta H_B' = (1 \otimes H_B' + H_B' \otimes G_B' F_B') \Delta
\end{equation}
The side conditions for $(F_{sym}, G_{sym}, H_{sym})$ and
$(F_B', G_B', H_B')$ follow from those for $(F, G_1, H)$.

\bigskip
\noindent
Next we  replace $(\mathcal{L}_1, L_1)$ by $(\mathcal{L}_2, L_2)$, 
taking into account  the Lie structures. 
The symmetric DG bialgebras of $\mathcal{L}_1$ and $L_1$ turn
into the universal enveloping DG bialgebras $U(\mathcal{L}_2)$, 
$U(L_2)$, respectively. Denote by $\rho: Sym^\bullet(\cdot) \to U(\cdot) $ the Poincare-Birkhoff-Witt 
isomorphism which identifies the two spaces as DG coalgebras, cf. 
Propositions 21.2 and 22.6 in \cite{FHT}.  Denote by $*$ the product in the 
universal enveloping and by $\cdot$ the product in the symmetric bialgebra. 

Since both $\mathcal{L}_2$ and $L_2$ are 2-step nilpotent (i.e. all double
brackets vanish) it is easy to track the multiplicative behavior of $\rho$.
In fact, let $K$  be a general graded Lie algebra with bracket $l_2$ such that
$l_2(l_2(a, b), c)=0$ for all $a, b, c \in K$. Then for 
odd elements $v_i, w_j \in K$
\begin{equation}\label{product}
\rho(v_1 \cdot \ldots \cdot v_n) *  \rho(w_1 \cdot \ldots \cdot w_m)
= \hspace{10cm}
\end{equation}
$$\sum_{k \geq 0} 
\sum_{\genfrac{}{}{0pt}{}
{I = \{i_1, \ldots, i_k\} \subset \{1, \ldots, n\}}
{J = \{j_1, \ldots, j_k\} \subset \{1, \ldots, m\}}}
(-1)^{(I, J)}  
det\Big\vert\frac{1}{2}l_2(v_{i_p}, w_{j_q})\Big\vert_{_{p, q = 1, \ldots, k }} *
\rho(v_{r_1} \cdot \ldots \cdot v_{r_{n-k}} \cdot w_{l_1}
\cdot \ldots \cdot w_{l_{m-k}})_{_{r_t \notin I, \;l_s \notin J}}
$$
where $I, J$ are subsets of equal cardinality with the
induced natural ordering and 
 $(-1)^{(I, J)}$ is defined by the equality in $Sym^\bullet(K)$:
$$
v_1 \cdot \ldots \cdot v_n \cdot w_1 \cdot \ldots \cdot w_m
= (-1)^{(I, J)}
v_{r_1} \cdot \ldots \cdot v_{r_{n-k}} \cdot v_{i_k} \cdot 
\ldots \cdot v_{i_1} \cdot w_{j_1} \cdot 
\ldots \cdot w_{j_k} \cdot w_{l_1} \cdot \ldots \cdot w_{l_{m-k}} 
$$
This formula is proved by first considering the case $n=1$ where 
it reduces to an easy computation, and then iterating and (anti)symmetrizing in $v_1, \ldots, v_n$.

Using the PBW isomorphism we can view $(F_{sym}, G_{sym}, H_{sym})$
as a contracting homotopy between $U(\mathcal{L}_2)$ and $U(L_2)$ but 
now $F_{sym}$ will not be multiplicative since $F: \mathcal{L}_2 \to L_2$
is not a Lie map. To ``repair" this we adjust the homotopy
on the bar construction.

Let $D_B, d_B$ be the canonical differentials of $BU(\mathcal{L}_2)$, $BU(L_2)$
and denote by $D_B', d_B'$ another pair of differentials on the same 
spaces, arising from their PBW isomorphism with $B Sym^\bullet(\mathcal{L}_2)$
and $B Sym^\bullet(L_2)$, respectively. The contracting homotopy $(F_B', G_B', H_B')$ 
may be viewed as a homotopy between $(BU(\mathcal{L}_2), D_B')$ and
$(BU(L_2), d_B')$.  Using the Basic Perturbation Lemma we can adjust it
to work with $D_B, d_B$ as follows.

Let $\delta_B = D_B - D_B'$. Explicitly, $\delta_B$ is obtained by 
considering the map
$$
\delta_U: U(\mathcal{L}_2) \otimes U(\mathcal{L}_2) \to U(\mathcal{L}_2), \qquad 
\rho(a) \otimes \rho(b) \mapsto \rho(a) * \rho(b) - \rho(a \cdot b)
$$
and then extending to $B U(\mathcal{L}_2)$  as a coderivation. Now set
$$
X = \delta_B - \delta_B H'_B \delta_B + \delta_B H'_B \delta_B H_B' \delta_B - \ldots.
$$
This infinite expression is well defined since $\delta_B$ decreases the tensor
degree in $BU(\mathcal{L}_2)$ by 1 and $H_B'$ preserves this degree. By 
Basic Perturbation Lemma, cf. e.g. \cite{GLS}, the formulas
$$
F_B = F_B'(1 - XH_B'); \quad H_B = H_B'(1 - X H_B'); \quad G_B = (1 - H_B' X) G_B'
$$
define a contracting homotopy between complexes $(BU(\mathcal{L}_2), D_B)$ and
$(BU(L_2), d_B' + F_B' X G_B')$
\begin{prop} The following properties hold
\begin{enumerate}
\item $\delta_B G_B' = G_B' (d_B - d_B')$, $G_B = G_B'$, $d_B = d_B' + F_B' X G_B'$;
\item $F_B$ is a coalgebra map and $H_B$ is a coalgebra homotopy,
see \eqref{ch};
\item $F_B$ and $H_B$ are uniquely determined by the compositions
$$
F_{B, k}: (U^{\geq 1} (\mathcal{L}_2))^{\otimes k} \to B U(\mathcal{L}_2) 
\to  BU(L_2) \to U^{\geq 1}(L_2) ;
$$
$$
H_{B, k}: (U^{\geq 1} (\mathcal{L}_2))^{\otimes k} \to B U(\mathcal{L}_2) 
\to  BU(\mathcal{L}_2) \to U^{\geq 1}(\mathcal{L}_2).
$$
\end{enumerate}
\end{prop}
\textit{Proof.} The identity $\delta_B G_B' = G_B' (d_B - d_B')$ follows 
from the fact that $G_1: L_2 \to \mathcal{L}_2$  commutes with brackets.
The other two identities follow from it and a side condition
$H_B' G_B' = 0$.

Part (2) is proved in \cite{GLS}. Part (3) is an easy consequence: for $F_B$
it is well-known, cf. e.g. \cite{Ka2}, while for $H_B$ one has 
an explicit formula
$$
H_B|_{(U^{\geq 1} (\mathcal{L}_2))^{\otimes k}} = \sum_{s + p + q = k} 1^{\otimes s}
\otimes H_{B,p} \otimes (G_B F_B)^{\otimes q}. \qquad \square
$$

\bigskip
\noindent
To summarize the above: we have defined a map of DG bialgebras $G_{sym}:
U(L_2) \to U(\mathcal{L}_2)$, an $A_\infty$-map of associative algebras
$\{F_{B, k}: U(\mathcal{L}_2)^{\otimes k} \to U(L_2)\}_{k \geq 1}$, and a system of 
``higher homotopies" $\{H_{B, k}: U(\mathcal{L}_2)^{\otimes k} \to U(\mathcal{L}_2)\}_{k \geq 1}$
which are encoded in a coalgebra contraction $(F_B, G_B, H_B)$ 
from $BU(\mathcal{L}_2)$ to $BU(L_2)$.

It is easy to see that $F_{B, k}$ and $H_{B, k}$ are given by (restrictions  of) 
$(-1)^{k-1} F_{sym} (\delta_B H'_B)^{k-1}$ and $(-1)^{k-1} H_{sym} (\delta_B H'_B)^{k-1}$,
respectively.

\bigskip
\noindent
Now we want to pass from $W^{(2)}$ to the full potential $W$. This means that the differential 
$D_B$ on $BU(\mathcal{L}_2)$ will be replaced by $\widehat{D}_B = D_B + \widetilde{\delta}_B$.
To describe $\widetilde{\delta}_B$ we first define a differential $\widetilde{\delta}_L$ on $\mathcal{L}_2$
which sends $v \in V \subset \mathcal{L}_2$ to $\widehat{\partial}_v (W - W^{(2)})$ (and vanishes
on the natural complement to $V$); then exend
$\widetilde{\delta}_L$ to a derivation $\widetilde{\delta}_U$ on $U(\mathcal{L}_2)$; and finally
extend $\widetilde{\delta}_U$ to a coderivation $\widetilde{\delta}_B$ on $BU(\mathcal{L}_2)$.
Using the Basic Perturbation Lemma again, we set
$$
\widetilde{X} = \widetilde{\delta}_B - \widetilde{\delta}_B H_B \widetilde{\delta}_B + 
\widetilde{\delta}_B H_B \widetilde{\delta}_B H_B \widetilde{\delta}_B - \ldots;
$$
which is well-defined since $\widetilde{\delta}_B$ decreases by 1 the number of occurences
of elements in $V \subset \mathcal{L}_2 \subset U(\mathcal{L}_2)$; and define
$$
\widetilde{F}_B = F_B(1 - \widetilde{X}H_B); \quad H_B = H_B(1 - \widetilde{X} H_B); \quad 
\widetilde{G}_B = (1 - H_B \widetilde{X}) G_B; \quad \widetilde{d}_B = d_B + F_B \widetilde{X}G_B.
$$
Then $(\widetilde{F}_B, \widetilde{G}_B, \widetilde{H}_B)$
is a contraction of $(BU(\mathcal{L}_2), \widetilde{D}_B)$ to $
(BU(L_2), \widetilde{d}_B)$. As in the previous Proposition we conclude that
$\widetilde{F}_B$ and $\widetilde{G}_B$ and maps of DG coalgebras, 
$\widetilde{H}_B$ is a coalgebra homotopy and $\widetilde{d}_B$ is a coderivation.

In particular, the coderivation $\widetilde{d}_B$ defines an $A_\infty$-structure on $U(L_2)$, cf. \cite{Le},
i.e a series of higher products 
$\mu_n: U(L_2)^{\otimes n} \to U(L_2)$ given by the 
composition of natural maps
$$
U(L_2)^{\otimes n} \to U^{\geq 1}(L_2)^{\otimes n}
\hookrightarrow BU(L_2) \stackrel{\widetilde{d}_B}\longrightarrow
BU(L_2) \to U^{\geq 1}(L_2) \hookrightarrow U(L_2). 
$$
Writing out the 
definitions and using $F_B' \widetilde{\delta}_B = 0$, $H_B' \delta_B G_B'  = 0$ 
we see that for $n \geq 3$, $\mu_n$ is given by the expression
\begin{equation}\label{higher}
 \sum_{k \geq (n-2); a_1, \ldots, a_k} (-1)^{k-1} F_{sym} \delta_U H_B' (a_1 H'_B)\ldots (a_k H_B') \widetilde{\delta}_B G_{sym}^{\otimes n} 
\end{equation}
where each $a_i$ is either  $\delta_B$ or $\widetilde{\delta}_B$ and the first possibility occurs 
precisely $(n-2)$ times. 

Alternatively, one can write a formula in the spirit of \cite{KS}:
$\mu_n$ is given by the sum over all planar trees with $n$ leaves,
one root and internal vertices of valency 2 or 3. Similarly
to \textit{loc. cit} we place $G_{sym}$ on each leaf, $F_{sym}$
on the root, $\delta_B$ on each internal vertex of valency 3,
$\widetilde{\delta}_B$ on each internal vertex of valency 2, and
$H_{sym}$ in the middle of each internal edge. A tree marked in 
this way is viewed as a ``flowchart" of operations applied to 
the arguments of $\mu_n$. Note that due to the valency 2
vertices each $\mu_n$ becomes an infinite sum 
over trees, but on each particular set of $n$ arguments only 
finitely many give nonzero contributions.

\begin{prop}  \label{properties}
The product $\mu_2$ is the usual unversal enveloping
algebra product in $U(L_2)$.
The higher products $\mu_n$ for $n \geq 3$ have the following
properties:

(1) Each $\mu_n$ is multilinear in $R = Sym^\bullet(U) \subset U(L_2)$.

(2) $\mu_n(a_1, \ldots, a_n) = 0$ if $a_i =1$ for some $i$. Thus, 
the $A_\infty$-structure is strictly unital.

(3) $\mu_n(v_1, \ldots, v_n) = \frac{1}{n!}l_n(v_1, \ldots, v_n)$
 if $v_i \in V \subset L \subset U(L)$ for all $i$. 
\end{prop}
\textit{Proof.}  To prove the assertion about $\mu_2$ first note that
the ``correction" to the product on $U(L_2)$ introduced by $\widetilde{d}_B
- d_B$ is given by the formula similar to 
as above expression for $\mu_n$, $n \geq 3$:
$$
\sum_{k \geq 1} (-1)^k F_{sym} \delta_U (H'_B \widetilde{\delta}_B)^k G_{sym}^{\otimes 2}
$$
but a single application of $H'_B \widetilde{\delta}_B$ will 
produce terms in $Sym^{\geq 1}(\{S\}) \subset U(\mathcal{L}_2)$. 
Since such terms are central in $U(\mathcal{L}_2)$, $\delta_U$ is
multilinear with respect to them.  But $F_{sym}$ vanishes on 
$Sym^{\geq 1}(\{S\})$ therefore the correction to the product on 
$U(L_2)$ vanishes.

Part (1) follows from \eqref{higher} (or better, the sum over
trees presentation) and the fact that the operators $F_{sym}$, 
$H_{sym}$ and $G_{sym}$ involved in it, are  all $R$-linear. Part
(2) follows from the fact that we are using the \textit{reduced} 
bar construction hence by definition all higher products factor 
through $U^{\geq 1}(L_2)^{\otimes n}$. To prove part (3) 
use the formula \eqref{higher} to compute 
$\mu(v_1, \ldots, v_n)$. 
The only non-zero contributions come 
from the terms with $k = (n - 2)$, i.e. for which all 
$a_i = \delta_B$. In fact, if a term in \eqref{higher} contains 
$\widetilde{\delta}_B$ at least twice, its evaluation at 
$v_1 \otimes \ldots v_n$ will necessarily contain 
$\delta_U (a \otimes b)$ with 
$a, b \in \{S\} \subset U(\mathcal{L}_2)$. To explain
that in terms of trees: 
if we connect the two occurences of $\widetilde{\delta}_B$
on a tree with its root by shortest paths, 
the point at which the two paths merge will correspond 
to the $\delta_U (a \otimes b)$ above. Since $\{S\}$ is 
central,  $\delta_U (a \otimes b) =0$. Therefore
$$
\mu_n (v_1, \ldots, v_n) = (-1)^{n-1} F_{sym} 
(\delta_B H_B')^{n-1} \widetilde{\delta}_B G_{sym}^{\otimes n}
(v_1 \otimes \ldots \otimes v_n)
$$
Now an easy induction involving \eqref{product} finishes the 
proof. $\square$

\bigskip
\noindent
\textbf{Remark.} The properties stated in the previous proposition do not
determine the $A_\infty$-structure uniquely. In the case of projective complete
intersections, cf. \cite{Ba}, an additional formula allows to compute all 
higher operations recursively. Such a formula can be proved in this case as well
but this will not be done here.

By a recent work of Merkulov, cf. \cite{Me}, the $L_\infty$-structure on 
$L$ deforms the commutative and cocommutative bialgebra structure of
$U(L_1)$ to a structure of a homotopy bialgebra. However, this structure
depends on the choice of a minimal model of the bialgebra PROP, and 
a certain lift of a morphism of PROPs (see \cite{Me} for more detail). 
 
Comparing our construction with the standard bialgebra $\Omega C(L)$, one can
show that in the situation of this paper the higher products on $U(L_2)$ extend to a homotopy bialgebra structure and in fact determine it uniquely.
We plan to return to this matter in a forthcoming work.

\bigskip
\noindent
\textbf{Notation.} From now on we denote by $E_W$ the universal enveloping $U(L_2)$ equipped with the $A_\infty$-algebra structure of this section.

\section{An equivalence of categories}

\subsection{A generalized twisted cochain.}

Let $C_W = Ker(\delta_C) \cap Sym^\bullet(V) \subset C$ be
the ``dual coalgebra" of the polynomial quotient $S_W$ defined
in Section 2.2.
Write $C_W = k \oplus \overline{C}$ where $\overline{C} = 
Ker (\delta_C) \cap Sym^{\geq 1}(V)$. If $\overline{\Delta}: 
\overline{C} \to \overline{C} \otimes \overline{C}$ is the
reduced coproduct $(\Delta - Id \otimes 1 - 1 \otimes Id)$, and $\overline{\Delta}^{(k)}: \overline{C} 
\to \overline{C}^{\otimes k}$ are its iterations, then
$\overline{C} = \cup Ker(\overline{\Delta}^{(k)})$, i.e.
$C_W$ is \textit{cocomplete}. In fact, this property holds for 
the free cocommutative coalgebra
$C$ and $C_W$ is its subcoalgebra.

Denote by $\Delta^{(k)}:C \to C^{\otimes k}$ similar iterations
for $k \geq 2$ and set $\Delta^{(1)}$ to identity.
Consider the composition 
$\tau_W: C_W \hookrightarrow C \to L \hookrightarrow U(L) = E_W$.
In other words, we compose the projection
$C_W \to C_W \cap V$ with the embedding
$V \subset L \subset E_W$.

\begin{lemma}\label{gentwist}
The map $\tau_W$ satisfies the generalized twisted cochain 
condition, cf. Section 4.1 of \cite{Le}, which reads
in our case:
$$
\sum_{s \geq 1} \mu_s \circ \tau_W^{\otimes s} \circ 
\Delta^{(s)} = 0
$$
\end{lemma}
\textit{Proof.} Note that the infinite sum is well defined
since $\tau_W|_k = 0$ and $C_W$ is cocomplete. First consider
$
C \to L\hookrightarrow E_W.
$
Then by the last part of Proposition \ref{properties}
one has
$$
\tau \circ\delta_C  +\sum_{s \geq 1} \mu_s 
\circ \tau^{\otimes s} \circ \Delta^{(s)} = 0.
$$
Since $C_W \hookrightarrow C$ is a morphism of coalgebras,
the assertion for $C_W$ follows trivially. $\square$

\subsection{A pair of adjoint functors.}

The previous lemma allows to apply the general 
formalizm outlined in Sections 2.2.1 and 4.3.1 of \cite{Le}.
Since some of the formulas are given in \cite{Le} only for
DG-algebras we give the definitions here for reader's 
convenience. See \cite{Le} for definitions and properties
of $A_\infty$-algebras and modules over them.
 
Consider a general cocomplete coaugmented DG-coalgebra 
$(C, \delta_C)$, a strictly unital  $A_\infty$-algebra $E$ 
with $\delta_E= \mu_1^E$ and a generalized twisted cochain 
$\tau: C \to E$ satisfying
$$
\delta_C \circ \tau + \tau \circ \delta_E +\sum_{s \geq 2} \mu_s 
\circ \tau^{\otimes s} \circ \Delta^{(s)} = 0 
$$
Let $(N, \delta_N)$ be 
a counital DG comodule over $C$
with the reduced coaction map $\overline{\Delta}_N: N \to N \otimes 
\overline{C}$. We assume that $N$ is also cocomplete, i.e. 
$N = \bigcup_{s \geq 1} Ker\overline{\Delta}_N^{(s)}$, where 
$\overline{\Delta}^{(s)}_N: N \to N \otimes
\overline{C}^{\otimes (s-1)}$ is the reduction of the iterated 
coaction map $\Delta^{(s)}_N:   N \to N \otimes C^{\otimes (s-1)}$.
Whenever we speak of a \textit{filtered} morphism of 
cocomplete comodules, we always have in mind the filtration
by $Ker\overline{\Delta}_N^{(s)}$.

\bigskip
\noindent
Denote by  $\mathcal{F}(N)$ the tensor product
$N \otimes E$ with the differential 
$$
\delta_{\mathcal{F}(N)} =
\delta_N \otimes 1 + 1 \otimes \delta_{E} + 
\sum_{s \geq 2} (1 \otimes \mu_s^E)(1 \otimes \tau^{\otimes (s-1)} \otimes 1) 
(\Delta_N^{(s)} \otimes 1)
$$
which is well-defined since $N$ is cocomplete and $E$ is strictly 
unital. Then $\delta_{\mathcal{F}(N)}^2 = 0$ 
by the generalized twisted cochain condition. Also, 
 $\mathcal{F}(N)$ is an $A_\infty$-module
over $E$ with the action maps
$$
\mu^{\mathcal{F}(N)}_k: \mathcal{F}(N) \otimes E^{\otimes (k-1)} \to \mathcal{F}(N);
\quad (n \otimes a) \otimes a_1 \otimes \ldots \otimes a_{k-1} \mapsto
n \otimes \mu_k^E(a, a_1, \ldots, a_{k-1}) 
$$
for $k \geq 2$. This module structure is \textit{strictly
unital}: $\mu^{\mathcal{F}(N)}_2(x, 1_{E}) = x$ and 
$\mu^{\mathcal{F}(N)}_k
(x, a_1, \ldots, a_{k-1}) = 0$ if $k \geq 3$ and $a_i = 1_{E}$ 
for some $i$, since the same property was assumed about $E$.
If $\psi: N_1 \to N_2$ is a morphism of $C_W$-comodules then
$\mathcal{F}(\psi) = \psi \otimes 1: N_1 \otimes E \to N_2 \otimes E$
is a strict morphism of $E$-modules (i.e. commutes with all higher
products).

\bigskip
\noindent
In the other direction, take a strictly unital $A_\infty$-module 
$(M, \delta_M, \mu^M_k)$ over $E$, where 
$\mu^M_k: M\otimes E^{\otimes(k-1)} \to M $
are the action maps for $k \geq 2$, and
consider the $C$-comodule  $\mathcal{G}(M) = M \otimes C$, 
 with the differential
$$
\delta_{\mathcal{G}(M)}= \delta_M \otimes 1 + 
1 \otimes \delta_C + \sum_{k \geq 2} (\mu^M_k  \otimes 1)
(1 \otimes \tau_W^{\otimes (k-1)} \otimes 1) 
(1 \otimes \Delta^{(k)}).
$$
Again the differential is well-defined since $C$ is cocomplete and
$E$ is strictly unital. A morphism of $A_\infty$-modules $M_1$,
$M_2$ is given by degree $(1-k)$ maps $f_k: M_1 \otimes E^{\otimes (k-1)}
\to M_2$ for $k \geq 1$,
which satisfy some quadratic identities, cf. Chapter 2 of \cite{Le}.
Such a morphism $f_{\cdot} = \{f_k\}$ is called \textit{strictly unital}
if $f_k(m, a_1, \ldots a_{k-1}) = 0$ whenever $k \geq 2$ and $a_i = 1_E$
for some $i$. For every such morphism define a morphism of 
$C$-comodules $\mathcal{G}(f_\cdot): M_1 \otimes C \to M_2 \otimes C$
by the formula
$$
\mathcal{G}(f_\cdot) = \sum_{k \geq 1} (f_k \otimes 1) 
(1 \otimes \tau^{\otimes (k-1)} \otimes 1)(1 \otimes \Delta^{(k)}). 
$$
This is well-defined for the same reason as before.

\bigskip
\noindent
Thus we obtain a pair of functors $\mathcal{F}$, $\mathcal{G}$ between
the category $Comodc(C)$ of cocomplete counital DG-comodules over $C$ and
the category $Mod_\infty(E)$ of strictly unital $A_\infty$-modules over 
$E_W$ and strictly unital morphisms. These functors are adjoint:
$$
Hom_{Mod_\infty(E)}(\mathcal{F}(N), M) = Hom_{Comodc(C)}(N, \mathcal{G}(M))
$$
since both spaces may be identified with
$$
\big\{\phi \in Hom_k (N, M) \;|\; \phi \delta_N - \delta_M \phi =
\sum_{k \geq 2} \mu_k^M (\phi \otimes \tau^{\otimes (k-1)}) \Delta^{(k)} \big\}.
$$
More explicitly, given such $\phi$ one defines a morphism
of $C$-comodules $\Phi: N 
\stackrel{\Delta_N}\longrightarrow N \otimes C 
\stackrel{\phi \otimes 1}\longrightarrow M \otimes C$,
and a morphism of $A_\infty$-modules 
$\Psi_\cdot: N \otimes E \to M$
$$
\Psi_k = \sum_{s \geq 1} \mu_{k+s}^M (\phi \otimes \tau^{\otimes(s-1)} 
\otimes 1^{\otimes k})(\Delta_N^{(s)} \otimes 1^{\otimes k}):
N \otimes E^{\otimes k} \to M.
$$
The map $\phi$ may be recovered from $\Psi_\cdot$ as 
$\Psi_1|_{N \otimes 1}$, or from $\Phi$ as its composition with 
the projection $\eta_M: M \otimes C \to M \otimes k = M$  coming
from the counit of $C$. The fact that the above formulas
indeed define morphisms, and that 
every $\Phi$, $\Psi_\cdot$
is given by a certain $\phi$, is proved by a straightforward
(but tedious) induction using the filtration of $N$ by 
$Ker(\overline{\Delta}^{(k)}:
N \to N \otimes \overline{C}^{\otimes (k-1)})$.

Below we need an explicit formula for the adjunction morphism
$\Psi_\cdot: \mathcal{F}\mathcal{G}(M) \to M$. The component
$\Psi_k: M \otimes C \otimes E^{\otimes k} \to M$ is given by 
$$
\sum_{s \geq 1}  \mu_{k+s}  (\eta_M \otimes \tau^{\otimes (s-1)} \otimes 1)
(1 \otimes \Delta^{(s)} \otimes 1)
$$

\subsection{A coalgebra equivalence.}

The generalized twisted cochain $\tau_W: C_W \to E_W$ extends to
a  coalgebra map $C_W \to B(E_W)$ by the standard formula 
$\sum_k \tau_W^{\otimes k} \Delta^{(k)}$, cf. \cite{Ka2}. 
The condition of Lemma \ref{gentwist} is equivalent to the fact
that this extension commutes with differentials. 

\begin{lemma} 
The canonical coalgebra extension $C_W \to B(E_W)$ defined by 
$\tau_W$, is a weak equivalence of coalgebras, i.e. induces a 
quasi-isomorpism of DG algebras $\Omega(C_W) 
\to \Omega B(E_W)$.
\end{lemma}
\textit{Proof.} Recall that the $A_\infty$-structure on $E_W$
is encoded in the differential $d_B$ on $B(E_W)$. In 
the previous section we have also constructed a quasi-isomorphism
of DG coalgebras $\widetilde{F}_B: BU(\mathcal{L}) \to BE_W$
which naturally induces a  
quasi-isomorphism of DG algebras 
$\Omega B U(\mathcal{L}) \to \Omega B (E_W)$.
It follows from the definitions 
that the algebra homomorphism $\Omega(C_W) 
\to \Omega B(E_W)$ factors as
$$
\Omega(C_W) \to \Omega(C) \to \Omega B U(\mathcal{L})
\to \Omega B(E_W)
$$
where the first and the last arrows are quasi-isomorphisms.
Therefore it suffices to check that 
$\Omega(C) \to \Omega B U(\mathcal{L})$
is a quasi-isomorphism. To that end, we note that the composition
$$
\Omega(C) \to \Omega B U(\mathcal{L}) \to U(\mathcal{L}),
$$
is a quasi-isomorphism by Lemma 3, and the second arrow 
is a quasi-isomorphism by a standard result in homotopical algebra
(see e.g. page 272 of \cite{FHT}). Therefore the first
arrow is also a quasi-isomorphism, which finishes the proof. 
$\square$

\bigskip
\noindent
Let $\mathcal{D}(E_W)$ be the localization $Mod_\infty(E_W)$ 
at quasi-isomorphisms.
To get a derived category $\mathcal{D}(C_W)$ we must 
localize $Comodc(C_W)$ at weak equivalences, i.e. such maps that induce 
quasi-isomorphism on cobar construction, cf. \cite{Le}. In 
general, a weak equivalence of comodules is a stronger condition than quasi-isomorphism.

\begin{corr}\label{coequiv}
The functors $\mathcal{F}$, $\mathcal{G}$ induce 
mutually inverse derived equivalences 
$\mathcal{D}(C_W) \simeq \mathcal{D}(E_W)$ 
between the derived category $\mathcal{D}(C_W)$ of 
cocomplete comodules over $C_W$ and the derived category 
$\mathcal{D}(E_W)$ of strictly unital $A_\infty$-modules over $E_W$. 
\end{corr}
\textit{Proof.} We  factorize $\mathcal{F}$ and $\mathcal{G}$  as follows
$$
\begin{array}{ccccc}
 & \stackrel{\mathcal{F}_1}\longrightarrow & & \stackrel{\mathcal{F}_0}\longrightarrow  &  \\
Comodc(C_W)  & & Comodc (B(E_W)) & & Mod_\infty(E_W) \\
 & \stackrel{\longleftarrow}{\mathcal{G}_1} & &  \stackrel{\longleftarrow}{\mathcal{G}_0}& 
\end{array}
$$
Here $\mathcal{F}_0$ and $\mathcal{G}_0$ are induced by the 
universal generalized twisted cochain 
$B(E_W) \to \overline{E}_W  \hookrightarrow E_W$,
where $\overline{E}_W$ is the kernel of the augmentation 
map. 
The functor $\mathcal{F}_1$ is given by corestriction 
(i.e. every $C_W$-comodule
is automatically a $B(E_W)$-comodule); and $\mathcal{G}_1$ by 
coinduction:
$$
\mathcal{G}_1(N) = Ker\big\{\Delta_N \otimes 1 - 
(1 \otimes \tau \otimes 1)(1 \otimes \Delta_{C_W}):
N \otimes C_W \to N \otimes B(E_W) \otimes C_W\big\}.
$$
It follows from definitions that
$\mathcal{F} = \mathcal{F}_0 \mathcal{F}_1$ and
$\mathcal{G} = \mathcal{G}_1 \mathcal{G}_0$.
Therefore, to prove that for any object $M$ of $Mod_\infty(E_W)$
the first component of the canonical $A_\infty$-morphism
$\mathcal{F}\mathcal{G} (M) \to M$ is a quasi-isomorphism, one needs to show
that

1) For any object $L$ of $Comodc(B(E_W))$ the canonical morphism
$\mathcal{F}_1 \mathcal{G}_1(L) \to L$ is a weak equivalence;

2) For any $M$ as above the canonical morphism 
$\mathcal{F}_0 \mathcal{G}_0 (M) \to M$ is a quasi-isomorphism;

3) $\mathcal{F}_0$ sends weak equivalences to quasi-isomorphisms.

\medskip
\noindent
Similarly, to prove that for any object $N$ of $Comodc(C_W)$ the 
canonical morphism $N \to \mathcal{G} \mathcal{F} (N)$ is a weak
equivalence, one needs to show that 

1') For any $N$ as above the canonical morphism $N \to \mathcal{G}_1 
\mathcal{F}_1(N)$ is a weak equivalence;

2') For any $L$ as above $L \to \mathcal{G}_0 \mathcal{F}_0(L)$ 
is a weak equialence;

3') $\mathcal{G}_1$ sends weak equivalences to weak equivalences. 

\bigskip
\noindent
In addition, to prove that the functors descend to derived categories
one needs to show that

4) $\mathcal{G}_0$ sends quasi-isomorphisms to weak equivalences;

5) $\mathcal{F}_1$ sends weak equivalences to weak equivalences.

\medskip
\noindent
The assertions 1), 1'), 3') and 5) follow from the previous Lemma, 
since the morphism 
$\Omega(C_W) \to \Omega B(E_W)$ gives restriction-induction
functors which descend to equivalences of derived categories; and 
$\Omega \mathcal{F}_1(N)$ is canonically isomorphic to the 
$\Omega(B(E_W))$-module induced from the $\Omega(C_W)$-module
$\Omega(N)$, while for $\Omega \mathcal{G}_1(N')$ there is a canonical 
quasi-isomorphism to the restriction of  $\Omega(N')$ from $\Omega(C')$ 
to $\Omega(C)$.

Statements 2) and 2') are proved by first setting formally $W =0$ 
where they both reduce to standard facts about the bar construction 
of an associative algebra, and then applying the Basic Perturbation 
Lemma to derive the case of general $W$ (recall, cf. Chapter 2 of 
\cite{Le}, that in 2') it
suffices to prove a filtered quasi-isomorphism).  

Finally, to prove 3) and 4) we first note that by defition
$\mathcal{G}_0$ sends
quasi-isomorphisms to filtered quasi-isomorphisms, which are automatically
weak-equivalences (the prooof is as in Lemma 1.3.2.2 of \cite{Le}).
To prove the assertion for $\mathcal{F}_0$ observe that the adjunction
morphism $B(E_W) \to B \Omega B(E_W)$ defines an $A_\infty$-morphism
$H_\cdot: E_W \to \Omega B(E_W)$. Therefore, every module $M$ over 
$\Omega B (E_W)$ becomes an $A_\infty$-module over $E_W$ if we
set $\mu^M_k: M \otimes E_W^{\otimes (k-1)} \to M$ to be the obvious
composition $M \otimes E_W^{\otimes (k-1)} \stackrel{1 \otimes H_{k-1}} 
\longrightarrow M \otimes \Omega B (E_W) \to M$. 

In particular, for any $L$ as in 1) the cobar construction $\Omega(L)$
is an $A_\infty$-module over $E_W$. Now the maps $1_L \otimes H_i$
for $i \geq 1$ give an $A_\infty$ morphism $\mathcal{F}_0(L) \to 
\Omega(L)$ of modules over $E_W$ which is easily seen to be a quasi-isomorphism
(since $H_\cdot$ is a quasi-isomorphism).
If $L \to L'$ is a weak equivalence in $Comodc(B(E_W))$ then
$\Omega(L) \to \Omega(L')$ is a quasi-isomorphism of modules
over $\Omega B (E_W)$ (and also over $E_W$), therefore $\mathcal{F}_0(L)
\to \mathcal{F}_0(L')$ is also a quasi-isomorphism.  $\square$

\subsection{The graded case.}

Now suppose that $A$ is an abelian group 
(not necessarily  torsion-free)
generated by elements $\alpha_1, \ldots, \alpha_n$ where 
$n = \dim_k V$. We fix a basis
$x_1, \ldots, x_n$ in $V^*$ and consider the $A$-grading 
$Sym^\bullet(V^*) = 
\oplus_{\alpha \in A} Sym^\bullet_\alpha(V^*)$ such that 
$\deg_A(x_i) = \alpha_i$. Assume that all $Sym^\bullet_\alpha(V^*)$ 
are finite dimensional. In this case, if $A_+$ denotes
the semigroup generated in $A$ by $\alpha_1, \ldots, \alpha_n$
then $A_+ \cap (-A_+) = \emptyset$ since otherwise there is a non-trivial 
monomial $x^I$ with $\deg_A(x^I) = 0$ and all its powers will 
satisfy the same condition too.

Assume further that polynomials $W_1, \ldots, W_m$ are $A$-homogeneous
of degrees $\beta_1, \ldots, \beta_m$
(see the next section for examples) and that neither of them 
has terms linear in $x_1, \ldots, x_n$ (this is mostly 
to simplify notation; if this condition is not satisfied, 
in the argument below one can either pass to a smaller polynomial 
quotient of  $Sym^\bullet(V^*)$ or adjust some definitions).
The last assumption implies that $L$ and $E_W$ have zero 
differentials.

Under these assumptions the quotient ring $S_W$ is also $A$-graded
with finite dimensional components, and the same holds
for $Sym^\bullet(V)$ (since the space $V$ will have the dual basis 
with $A$-degrees $-\alpha_1, \ldots, -\alpha_n$). The action of
$Sym^\bullet(V^*)$ by differential operators on $Sym^\bullet(V)$ 
agrees with the $A$-grading and by assumption on $W_1, \ldots, W_m$
the coalgebra $C_W$ will inherit
the $A$-grading as well. It follows immediately from the 
definitions that the pairing 
$\langle\;,\;\rangle: Sym^\bullet(V^*) \times
Sym^\bullet(V) \to k$ descends to $\langle\;,\;\rangle: S_W \times
C_W \to k$.

The $A_\infty$-algebra $E_W \simeq \Lambda^\bullet(V) \otimes 
Sym^\bullet(U)$ has $A$-grading in which $z_1, \ldots, z_m$ have degrees
$-\beta_1, \ldots, -\beta_m$, respectively, (this ensures that
$\deg_A \mu_k = 0$ for $k \geq 1$). Assuming that the upper indices
refer to homological grading and lower indices to $A$-grading we see
that $(E_W)^i_\alpha = 0$ unless $(\alpha, i)$ is in
$$
\mathcal{A} = \{\textrm{subsemigroup of } A \times \mathbb{Z}
\textrm{ generated by } (-A_+, 1), (-A_+, 0) 
\textrm{ and } 0\}.
$$ 
Here we use the assumption that $W_1, \ldots, W_m$ have no linear terms. 

Let $C^\clubsuit(S_W)$
be the category of complexes $N = (\ldots \to N^{-1} \to N^0 \to N^1 \to \ldots)$ 
of $A$-graded modules $N^i = \oplus_{\alpha \in A} N^i_\alpha$ 
over $S_W$, such that $N_\alpha^i = 0$ unless 
$(\alpha, i) \in \alpha + \mathcal{A}$ 
for some $\alpha \in A \times \mathbb{Z}$ depending on $N$. Note that
objects in the category $C^b(S_W)$ of finite complexes of finitely 
generated $S_W$-modules will not in general satisfy this condition.
However, the \textit{opposite} category $(C^b(S_W))^{opp}$
(with inverted arrows) can embedded into $C^\clubsuit(S_W)$ if for any $N 
\in Ob(C^b(S_W))$ we consider the graded dual $N^*$ with 
$
\big(N^*\big)^{-i}_{-\alpha} = (N^i_\alpha)^*.
$

\bigskip\noindent
Similarly we define $C^\clubsuit(E_W)$ as the category of strictly unital
$A$-graded $A_\infty$-modules $M$ over $E_W$ which satisfy 
$M^i_\alpha = 0$
unless $(\alpha, i) \in \gamma + \mathcal{A}$ for some 
$\gamma \in A \times \mathbb{Z}$
depending on $M$.

\bigskip
\noindent
Any $S_W$-module $N \in Ob(C^\clubsuit(S_W))$ 
is also a $C_W$-comodule
with the reduced coaction map
$$
\overline{\Delta_N}: N \to N \otimes \overline{C}_W; \qquad n \mapsto
\sum_{\alpha \in A_+, y_{\alpha, i}} (n y_{\alpha, i}) 
\otimes x_{\alpha, i}
$$
where the sum is taken over dual bases 
$\{x_{\alpha, i}\}_{1 \leq i \leq \dim (C_W)_\alpha}$
and $\{y_{\alpha, i}\}_{1 \leq i \leq \dim (S_W)_\alpha}$ of 
$(C_W)_\alpha$ and $(S_W)_\alpha$, 
respectively. The condition imposed on the grading of $N$, together with 
$(-A_+) \cap A_+ = \emptyset$, ensure that the sum in the defintion 
of $\overline{\Delta}_N$ is finite on every $n \in N$.

Thus we can still define a functor 
$\mathcal{F}: C^\clubsuit(S_W) 
\to C^\clubsuit(E_W)$ sending $N$ to $N \otimes E_W$.
As for $\mathcal{G}: C^\clubsuit(E_W) 
\to C^\clubsuit(S_W)$, note that 
$\mathcal{G}(M) = M \otimes C_W$ is not only a 
$C_W$-comodule but also an $S_W$-module (since $C_W$ itself is 
a graded dual to the free rank one module over $S_W$). These 
functors descend to the corresponding derived categories 
$D^\clubsuit(S_W)$, $D^\clubsuit(E_W)$.

To formulate the next theorem we define the ``bounded" derived
category $D^b(E_W)$ as the full triangulated subcategory 
of $D^\clubsuit(E_W)$ formed by all objects for which the 
total cohomology (= direct sum of all cohomology groups) is a 
finitely generated module over the associative algebra 
$(E_W, \mu_2)$. Note that the $E_W$ itself is only bounded
from the left, so many objects in $D^b(E_W)$ will be 
unbounded in the usual sense. Let also $D^b(S_W)$ be 
the usual bounded derived category of finitely
generated $A$-graded  $S_W$-modules, embedded contravariantly
 into 
$D^\clubsuit(S_W)$ by the above.

\begin{theorem}\label{bgg} 
The functors $\mathcal{F}$, $\mathcal{G}$ give 
mutually inverse derived equivalences
$$
D^\clubsuit (S_W) \to D^\clubsuit(E_W)
$$
Moreover, their restrictions induce a derived equivalence
$$
D^b(S_W)^{opp} \simeq D^b(E_W)
$$
\end{theorem}
\textit{Proof.} From the proof of Corollary \ref{coequiv} we already know that 
the adjunction morphisms $\mathcal{F}\mathcal{G}(M) \to M$ and 
$N \to \mathcal{G}\mathcal{F}(N)$ are quasi-isomorphisms. To establish 
the first claim it remains to show that $\mathcal{F}$ and $\mathcal{G}$
send quasi-isomorphisms to quasi-isomorphisms. This is obvious for
$\mathcal{G}$ since it sends quasi-isomorphisms to weak equivalences
of $C_W$-comodules (again by proof of Corollary \ref{coequiv}) and
every weak equivalence is a quasi-isomorphism. To prove the
assertion for $\mathcal{F}$ first assume that $N \to N'$ is a
quasi-isomorphism in $C^\clubsuit(S_W)$. Then by the above
$\mathcal{G}\mathcal{F}(N) \to \mathcal{G}\mathcal{F}(N')$ is 
also a quasi-isomorphism. It follows from the definitions that 
the last map can be viewed as filtered quasi-isomorphism (hence
a weak equivalence) of $C_W$-comodules. Apply the proof of 
Corollary \ref{coequiv} again we see that $\mathcal{F} \mathcal{G}
\mathcal{F}(N) \to \mathcal{F} \mathcal{G} \mathcal{F}(N')$ is
a quasi-isomorphism and therefore $\mathcal{F}(N) \to \mathcal{F}(N')$
is a quasi-isomorphism. 

The proof of the second assertion proceeds
exactly as  in 
of Proposition 4.1 in \cite{Ba}: first we note that
for an object $N$ in $D^b(S_W)$ the cohomology of 
$\mathcal{F}(N)$ is simply $Ext^\bullet(N, k)$ therefore
by Section 3 of \cite{G} it is finitely generated over
$(E_W, \mu_2) \simeq Ext_{S_W}(k, k)$. Therefore, 
$\mathcal{F}$ sends $D^b(S_W)^{opp}$ to $D^b(E_W)$. 
In the other direction, if $M$ is an object of $D^b(E_W)$
we can first replace it by $\mathcal{F} \mathcal{G}(M)$ which 
is a complex of free $E_W$-modules. It suffices to
check that $\mathcal{G}\mathcal{F}\mathcal{G}(M)$ gives
a finitely generated module \textit{over $S = Sym^\bullet(V^*)$}
but that module may be computed using the equivalence
functors $\mathcal{F}_0$ and $\mathcal{G}_0$ for
$S$ and $E = \Lambda^\bullet(V)$, respectively. Note that
$E$ is a quotient of $E_W$(the quotient map sends all 
$z_j \in U$ to zero). Inspecting the definitions of the
functors involved we see that 
$\mathcal{G}\mathcal{F}\mathcal{G}(M)$ is isomorphic to 
$\mathcal{G}_0 \big(\mathcal{F}\mathcal{G}(M)\otimes_{E_W} E)$,
as objects of $D^b(S)^{opp}$. Therefore the finite generation
of the graded dual follows from the original BGG correspondence
for $S$ and $E$ (adjusted to the case of $A$-graded modules).
$\square$ 

\bigskip
\noindent
\textbf{Remark.} When $A = \mathbb{Z}$ and 
all $W_i$ are quadratic one can adjust the grading to ensure 
that $E_W$ is in homological degree zero. 
In the general case this will not be possible since by 
Proposition 4 the 
higher products $\mu_k$, $k \geq 3$ on $E_W$ will be nontrivial 
and they have homological degree $2-k$.

\section{Toric complete intersections and Landau-Ginzburg models.}

\subsection{Homogeneous coordinates on toric varieties}

We fix notation by recalling some facts about toric varities.
Let $T \simeq (\k^*)^n$ be  
 an algebraic torus over $\k$.  Associated to $T$ are
the character lattice $M = Hom_{alg} (T, \k^*)$ and the dual
lattice of one-parametric subgroups $N = Hom_{alg} (\k^*, T)$ (the 
subscript \textit{alg} means homorphisms 
in the category of algebraic groups over $k$). Set also
 $N_\R = N \otimes_\Z \R$
and consider a fan $\Sigma \subset N_\R$ defining a toric variety 
$X_\Sigma$, cf. \cite{F}. We assume that the support of $\Sigma$
is equal to $N_\R$, i.e. the variety $X_\Sigma$ is complete.
If $\Sigma(1)$ is the set of all 1-dimensional
cones in $\Sigma$ let $S$ be the polynomial algebra over $\k$
generated by variables $x_\rho$, $\rho \in \Sigma(1)$.

Recall, cf. \cite{F}, that $\Sigma(1)$ is in bijective correspondence
with the codimension 1 orbits of the $T$-action on $X_\Sigma$.
For any $\rho \in \Sigma(1)$ let $D_\rho$ be the closure of the
corresponding orbit. Define the group $A$ by
the exact sequence 
$$
0 \to M  \stackrel{\alpha}\to \bigoplus_{\rho \in \Sigma(1)}
\Z D_\rho \stackrel{\beta} \to A \to 0
$$
where $\alpha(m) = \sum_\rho \langle m , n_\rho \rangle D_\rho$, 
$\beta$ is the quotient map and $n_\rho \in N$ is the 
primitive generator of $\rho \subset N_\R$.  
If  we view the elements of $M$  as rational functions on $X_\Sigma
\supset T$ then $\alpha$ computes the orders of poles and 
zeros along the divisors
$D_\rho$. One can show that $A$ is isomorphic to the Chow group 
$A_{n-1}(X_\Sigma)$. In general $A$ will have torsion.

Note that the map $\beta$ gives an $A$-grading on $S$. For
projective spaces this reduces to the usual 
$\Z$-grading on polynomials. Denoting
$$
G = Hom_Z (A, \k^*)
$$
we get a $G$-action on the space $V = \k^{\Sigma(1)}$ dual 
to the vector space spanned by $x_\rho$, $\rho \in \Sigma(1)$. 
To obtain $X_\Sigma$
from this action first denote by $\Sigma(max)$ the set of
\textit{maximal} cones in $\Sigma$ (i.e. those which are not
contained in a larger cone) and then for $\sigma \in \Sigma(max)$
define $x^{\s{\sigma}} \in S$ as $\prod_{\rho \in \Sigma(1) \setminus
\sigma(1)} x_\rho$. The monomials $x^{\s{\sigma}}$, $\sigma \in \Sigma(max)$
generate an ideal $B \subset S$ which corresponds to a closed subvariety
$V(B) \subset V$. To describe $V(B)$ more explicitly,
let $\Gamma \subset \Sigma(1)$ be a subset and $V(\Gamma)$ the
coordinate subspace of $V$ defined by vanishing
of the coordinates in $\Gamma$. Then $V(\Gamma) \subset V(B)$
iff $\Gamma$ is not contained in the closure of any maximal 
cone $\sigma \in \Sigma(max)$. It is easy to see that $V(B)$
is the union of all such $V(\Gamma)$ and we can restrict this
union to those $\Gamma$ which are minimal (with respect to 
inclusion) among all subsets with the above property. 

Since the support of $\Sigma$ is $N_{\mathbb{R}}$, every 
$\Gamma$ which is not in the closure of a maximal cone must contain 
at least 2 elements. Therefore $codim_V V(B) \geq 2$. 

\bigskip
\noindent
According to the main result of \cite{C1} we have
$$
X_\Sigma = \big[V \setminus V(B) \big]/G
$$
The right hand side is usually understood  as the universal 
categorical quotient or geometric quotient, 
if $\Sigma$ is simiplicial
(i.e. all cones in $\Sigma$ are simplices). However, in this
paper we will view it as a \textit{stack}
and use the above equality to \textit{define} the quotient
stack $\X$. Thus, in general $\X$ will be an Artin stack;
when $\Sigma$ is simplicial, a Deligne-Mumford stack, cf. \cite{LM},
\cite{F}. 

\subsection{BGG correspondence for toric complete intersections}

Since $\X$ was explicitly defined as a global quotients
it is easy to deal with bundles and sheaves on $\X$: these are just
$G$-equivariant objects on $V^\circ = V \setminus V(B)$. 
For instance, line bundles on $\X$ are just $G$-equivariant line 
bundles on $V^\circ$. 

Since $V(B)$ has codimension 
at least two, $Pic(V^\circ)$ is trivial and $H^0(V^\circ, \mathcal{O})
= k$. It follows that the Picard group of $\X$ may be identified with 
$Hom_{alg} (G, k^*) = A$. Note that for
$X_\Sigma$ the Picard group is in general only a subgroup of $A$
(some line bundles on $\X$ give only torsion-free sheaves
on $X_\Sigma$). Thus,
for any $\alpha \in A$ we have a line bundle $\mathcal{O}(\alpha)$
on $\X$. If $S_\alpha \subset S$ is the graded component corresponding 
to $\alpha$ then 
$$
H^0(\X, \mathcal{O}(\alpha)) = S_\alpha
$$
and $\dim_k S_\alpha < \infty$, cf. \cite{C1}.

Now let $W_1, \ldots, W_m$ be a regular sequence of elements
in $S$ which are $A$-homogeneous of $A$-degrees $\beta_1, 
\ldots, \beta_m$, respectively. Since these can be viewed as 
sections of line bundles $\mathcal{O}(\beta_1)$, $\ldots$, 
$\mathcal{O}(\beta_m)$, respectively, they define a complete
intersection substack $\mathbb{Y}\subset \X$. In other words,
$$
\mathbb{Y} = [V^\circ \cap Z(J)]/G,
$$
where $J$ stands for the ideal generated in $S$ by 
$W_1, \ldots, W_m$ and $Z(J)$ is the zero set.  
Denote $S/J$ by $S_W$ (note that the notation for 
$J, S_W$ and $\alpha_i$, $\beta_j$ is consistent with 
that of Sections 2 and 3, respectively).

Recall that a full triangulated subcategory $I$ of a triangulated 
category  $D$ is called \textit{thick} if it is closed with 
respect to  the operation of taking direct summands. 

\begin{lemma}
The category $Coh(\X)$ of coherent sheaves on $\X$ is equivalent
to the quotient of the category $mod_A(S)$ of $A$-graded
finitely generated modules over $S$, by the subcategory of 
modules supported on $V(B)$. 

Similarly, the category $Coh(\mathbb{Y})$ of coherent sheaves 
on $\mathbb{Y}$ is equivalent the quotient of the category 
$mod_A(S_W)$ of $A$-graded finitely generated modules over
$C_W$, by the subcategory of modules supported on $V(B)\cap Z(J)$.
The derived category $D^b(Coh(\mathbb{Y}))$ is equivalent
to the quotient $D^b(S_W)/I$ where $I$ is the thick 
subcategory of all complexes with cohomology supported 
on $V(B)\cap Z(J)$. 
\end{lemma}
\textit{Proof.} For $Coh(\X)$,$Coh(\mathbb{Y}$ everything 
follows easily from definitions since a coherent
sheaf on $\X$, resp. $\mathbb{Y}$, is simply a $G$-equivariant  
coherent sheaf on $V^\circ$, resp. $V^\circ\cap Z(J)$
which can always be extended to $V$, resp $Z(J)$. 
Thus, $Coh(\mathbb{Y})$ is a quotient of $mod_A(S/J)$ and
the kernel is easily seen to  be the subcategory of modules
supported on $V(B)\cap Z(J)$. 
The derived category statement is similar. 
\hfill $\square$

\bigskip
\noindent
 Since the graded components of $S$ are 
finite-dimensional, by Theorem \ref{bgg}, the derived 
category $D^b(\mathbb{Y})$ is equivalent to a quotient
of the $D^b(E_W)^{opp}$. 
To describe this quotient, suppose that $\Gamma \subset \Sigma(1)$
defines an irreducible component of $V(B)$, i.e. that 
$\Gamma$ is a minimal subset of $\Sigma(1)$ which is not 
contained in the closure of a maximal cone of $\Sigma$.  
Set $L_\Gamma$ to be the subspace of $V^*$ spanned by 
$x_\sigma$ with $\sigma  \in \Gamma$. Let also $V_\Gamma$ be the 
annihilator of $L_\Gamma$ in $V$. Obviously, $L_\Gamma^* \simeq
V/V_\Gamma$. 

Recalling  the quotient $E_W \to 
E = \Lambda^\bullet(V)$ from the proof of Theorem \ref{bgg} 
we see that every 
$E$-module automatically becomes an $E_W$-module (with
vanishing higher products).

\begin{prop}
The derived category $D^b(\mathbb{Y})$ is equivalent to 
the $(D^b(E_W)/T)^{opp}$ where $T$ is the thick subcategory 
generated by the $A$-shifts of 
$E_W$-modules $\Lambda^\bullet(V/V_\Gamma)$ and $\Gamma$ runs 
through the collection of all minimal subsets in $\Sigma(1)$ 
which are not contained in the closure of 
any maximal cone of $\Sigma$.
\end{prop}
\textit{Proof.} Fix a $\Gamma$ as above and consider
the Koszul complex $(\Lambda^\bullet (L_\Gamma) \otimes S_W, 
d_{Kos})$ of the (not necessarily regular) sequence of 
elements in $S_W$ given by the images of $x_\sigma$, $\sigma 
\in \Gamma$. Since it is exact on 
$Z(J) \setminus (Z(J)\cap V(B))$ and its zero 
cohomology is the algebra of functions on the scheme 
intersection $Z_J\cap V(B)$, by Lemma 1.2 in \cite{N} 
the thick subcategory of $D^b(S_W)$ generated by the $A$-shifts of 
$(\Lambda^\bullet (L_\Gamma) \otimes S_W, 
d_{Kos})$ is precisely the category formed by objects
with cohomology supported on $Z(J)\cap V(B)$. 

By Theorem 1.5 in \cite{N} the thick subcategory $I$ of 
the previous lemma is generated by the $A$-shifts of
$(\Lambda^\bullet (L_\Gamma) \otimes S_W, d_{Kos})$
for all $\Gamma$ as in the statement of the proposition.

It remains to prove that the functor $\mathcal{F}$ of Theorem 
\ref{bgg} takes $(\Lambda^\bullet (L_\Gamma) \otimes S_W, 
d_{Kos})$ to the $E_W$-module $\Lambda(V/V_\Gamma)$.  
By the same Theorem \ref{bgg} it suffices to show instead 
that $\mathcal{G}$ takes $\Lambda(V/V_\Gamma)$ to 
$(\Lambda^\bullet (L_\Gamma) \otimes S_W, 
d_{Kos})$ but that follows from the definitions. $\square$

\subsection{Categories of singularities and BGG correspondence}

Let $Z= Z(J)$ be the complete intersection in $V$ as before
and $D^{perf}(Z)\subset D^b(Z)$ the triangulated subcategory of perfect 
complexes (in this case, complexes quasi-isomorphic
to finite complexes of finitely generated $A$-graded projective 
$S_W$-modules). The quotient $D_{sg}(Z) = D^b(Z)/D^{perf}(Z)$
 has been studied in \cite{O1}
and \cite{O2} in relation to Landau-Ginzburg models. 
See \textit{loc. cit.} for more details and motivation.
The following proposition is a direct consequence of 
Theorem \ref{bgg}.
\begin{prop}
The $A$-graded category $D_{sg}(Z)$ of singularities on $Z$
is equivalent to the quotient $D^b(E_W)/R$ where $R$ is 
the thick subcategory generated by the $A$-shifts of $k$. 
\end{prop}
\textit{Proof.} Recall, that the functor $\mathcal{F}$
simply sends an $S_W$-module $M$ to $Ext^\bullet(M, k)$
(viewed as an $A_\infty$-module over the Yoneda
algebra $Ext^\bullet(k, k)$). Therefore, if $M$ is projective,
$\mathcal{F}$ is a direct sum of several copies of $k$ with 
their $A$-grading shifted. 
The assertion folows. $\square$

\subsection{The Calabi-Yau case}

Suppose that $\mathbb{Y}$ has trivial canonical class. By
the canonical class formula in Section 4.3 of \cite{F} 
combined with the adjunction formula, this 
condition is equivalent to 
$$
\beta_1 + \ldots + \beta_m = \alpha_1 + \ldots + \alpha_n
$$ 
Suppose further that the set-theoretic intersection $Z(J)
\cap V(B)$ consists only of the origin. This is automatically 
the case when $\X$ is a weighted projective space (since then 
$V(B)$ itself reduces to the origin). In this case we have
an alternative description of $D^b(\mathbb{Y})$ as $D_{sg}(S_W)$,
cf. Theorem 2.5 in \cite{O2}. 
\begin{corr}
If $\mathbb{Y}$ has trivial canonical class and $Z(J) \cap V(B)$
is supported at the origin, there exists a derived equivalence
$$
D^b(\mathbb{Y}) \simeq D^b(E_W)/R
$$
where $R$ is the thick subcategory generated by $k$. $\square$
\end{corr}
When $\mathbb{Y}$ is a complete intersection of quadrics in a projective
space, after a slight adjustment of grading on $E_W$ 
(see the Remark at the end of Section 3) this 
reduces to the result of Bondal and Orlov, \cite{BO}.

\noindent
\textsl{Department of Mathematics, 103 MSTB}\\
\noindent
\textsl{University of California, Irvine}\\
\noindent
\textsl{Irvine, CA 92697, USA}\\
\noindent
\textsl{email: vbaranov@math.uci.edu}


\begin{thebibliography}{BKRS}



\bibitem[Ba]{Ba} Baranovsky V.: BGG correspondence for projective complete
intersections, \textit{Internat. Math. Res. Notices}, \textbf{45} (2005), p.
2759-2774.

\bibitem[BGG]{BGG} Bernstein, I. N.; Gelfand, I. M.; Gelfand, S. I.:
Algebraic vector bundles on $ P^{n}$ and problems of linear algebra. (Russian)  
\textit{Funktsional. Anal. i Prilozhen.} \textbf{12}  (1978), no. 3, 66--67, 
English translation in  \textit{Functional Anal. Appl.} \textbf{12} (1978), 
no. 3, 212--214.

\bibitem[BGS]{BGS}
Beilinson, A.; Ginsburg, V.; Schechtman, V.: Koszul 
duality.  \textit{J. Geom. Phys.} \textbf{5}  (1988),  no. 3, 
317--350.

\bibitem[BO]{BO} Bondal, A.; Orlov, D.: 
Derived categories of coherent sheaves, preprint
math.AG/0206295


\bibitem[C]{C1} Cox, D. : The homogeneous coordinate ring of a 
toric variety, \textit{J. Algebraic Geom.} \textbf{4} (1995), no. 1, 
17--50. 

\bibitem[F]{F} Fulton, W.: Introduction to toric varieties,
Annals of Mathematics Studies, 131. The William H. Roever 
Lectures in Geometry. \textit{Princeton University Press, 
Princeton, NJ}, 1993.

\bibitem[FHT]{FHT}F\'elix, Y.; Halperin, S.; Thomas, J.-C.:
Rational homotopy theory. Graduate Texts in Mathematics, 205. 
\textit{Springer-Verlag, New York}, 2001.

\bibitem[G]{G}
Gulliksen, T. H.: 
A change of ring theorem with applications to Poincar\'e series 
and intersection multiplicity.  \textit{Math. Scand.}  
\textbf{34}  (1974), 167--183. 

\bibitem[GJ]{GJ} Getzler, E.; Jones, J. D. S.: Operads, homotopy 
algebra and iterated integrals for double loop spaces,
preprint hep-th/9403055. 

\bibitem[GLS]{GLS}
Gugenheim, V. K. A. M.; Lambe, L. A.; Stasheff, J. D.:
 Perturbation theory in differential homological algebra. II.  \textit{Illinois J. Math.}
\textbf{35}  (1991),  no. 3, 357--373.


\bibitem[Ka1]{Ka} Kapranov, M.: On the derived category and $K$-functor of 
coherent sheaves on intersections of quadrics. (Russian)  \textit{Izv. Akad. 
Nauk SSSR Ser. Mat.}
\textbf{52}  (1988),  no. 1, 186--199;  English translation in  
\textit{Math. USSR-Izv.} \textbf{32}  (1989),  no. 1, 191--204. 

\bibitem[Ka2]{Ka2} Kadeishvili, T. On the cobar construction of a bialgebra.  
\textit{Homology Homotopy Appl.} \textbf{7}  (2005),  no. 2, 109--122.

\bibitem[KS]{KS} Kontsevich, M.; Soibelman, Y.:
Homological mirror symmetry and torus fibrations.  
\textit{Symplectic geometry and mirror symmetry (Seoul, 2000)},  
203--263, \textit{World Sci. Publ., River Edge, NJ}, 2001; 
also preprint math.SG/0011041.


\bibitem[Le]{Le} Lef\`evre-Hasegawa, K.: Sur les A-infini 
cat\'egories, preprint math.CT/0310337.

\bibitem[LM1]{LMa} Lada, T.; Markl, M.:  Strongly homotopy Lie 
algebras, \textit{Comm. Algebra}  \textbf{23}  (1995),  no. 6, 
2147--2161; also preprint hep-th/9406095.

\bibitem[LM2]{LM} 
Laumon, G.; Moret-Bailly, L.:
Champs alg\'ebriques.
Ergebnisse der Mathematik und ihrer Grenzgebiete. 3. Folge. 
A Series of Modern Surveys in Mathematics, \textbf{39}. 
\textit{Springer-Verlag, Berlin}, 2000.

\bibitem[Me]{Me} Merkulov S.A.: Quantization of strongly homotopy Lie 
bialgebras, preprint math.AG/0612431.

%\bibitem[Ma]{Ma} Markl, M.:

\bibitem[N]{N}
Neeman, A.:
The chromatic tower for $D(R)$.
With an appendix by Marcel B\"okstedt.
\textit{Topology} \textbf{31} (1992), no. 3, 519--532. 


\bibitem[O1]{O1} Orlov, D.: Triangulated categories of singularities and equivalences between Landau-Ginzburg models, preprint
math.AG/0503630.

\bibitem[O2]{O2} Orlov, D.: Derived categories of coherent 
sheaves and triangulated categories of singularities, preprint
math.AG/0503632. 



\end{thebibliography}
\end{document}